\documentclass[11pt,twoside,reqno]{amsart}
 \usepackage{pgfplots}
\usepackage[colorlinks=true,citecolor=blue]{hyperref}
 \usepackage{mathptmx, amsmath, amssymb, amsfonts, amsthm, mathptmx, enumerate, color}
 \usepackage{graphicx}
 \usepackage{epstopdf}
 \newtheorem{theorem}{Theorem}[section]
 
 \newtheorem{lemma}{Lemma}[section]

 \newtheorem{remark}{Remark}[section]

 \numberwithin{equation}{section}

\begin{document}
 	\setcounter{page}{1}

 	\vspace*{1.0cm}
 	\title[A bilevel optimization  based algorithm  for solving\\ a class of price  equilibrium models ]{ A bilevel optimization  based algorithm
for solving\\ a class of price  equilibrium models \\}
 	 	\author[  Nguyen Ngoc Hai, Le Dung Muu, Nguyen Van Quy]{Nguyen Ngoc Hai$^1$, Le Dung Muu$^2$,  Nguyen Van Quy$^3$}
 	\maketitle
 	\vspace*{-0.6cm}

 	\begin{center}
 		{\footnotesize {\it
 				
                 $^1$ Trade Union University, Hanoi, Vietnam\\
 				$^2$TIMAS, Thang Long University, Hanoi, Vietnam and Institute of Mathmatics, VAST\\
                           
                            $^3$ Financial University, Hanoi, Vietnam}}

 	{\it  This paper dedicated to Professor Pham Ky Anh  on the occasion of his 75th birthday  }
 \end{center}
 	
 	\vskip 4mm {\small \noindent {\bf Abstract.} We consider class of equilibrium models including the implicit Walras supply-demand  and competitive  models.  Such a  model in this class, in general, is ill-posed.
We formulate such a model in the form a  variational inequality having certain monotonicity property   which allow us to describe a regularization algorithm avoiding the ill-posedness  based upon the bilevel optimization for finding a point that is nearest  to the given guessed or desired equilibrium price for the model.
  The obtained computational results   with many randomly generated data show that the proposed algorithm works well for this class of the equilibrium models.

 		   }

 	\renewcommand{\thefootnote}{}
 	\footnotetext{ $^1$Corresponding author.
 		\par
 		E-mail addresses: hainn@dhcd.edu.vn (N.N. Hai), ldmuu@math.ac.vn (L.D. Muu),  nguyenvanquy@hvtc.edu.vn (N.V. Quy).
        
        This paper has been submitted to the JORC on 24 Apr 2024.}
 	
 \newcommand{\N}{\mathbb{N}}

\section{Introduction} We consider  equilibrium models, where there are two sectors
with $n$-commodities $(x_1,...,x_n) \in \mathbb{R}^n$ depending on a parameter $p$ lying in a closed convex set $P$. Examples for such  models are 
the   Walras \cite{Wa1} or competitive equilibrium price ones. As usual, we suppose that each sector has a strategy set $S(p) \subset X \subseteq \mathbb{R}^n_+$ and $D(p) \subset X \subseteq \mathbb{R}^n_+$, and we call a point $p^*$  an equilibrium price of the model if there is $s^* \in S(p^*)$, $d^* \in D(p^*)$ such that  $\langle s^*-d^*, p-p^*\rangle \geq  0 \ \ \forall p\in P$, which means that $p^*$ is a solution of a multivalued variational  inequality problem.  In the case when $P$ is a closed convex cone,  this problem is reduced into a multivalued complementarity one. Complementarity problems even in the linear case, in general, is difficult to solve.
 
In this paper we suppose that  both strategies $S$ and $D$ are given implicitly as the optimal solutions  of certain parametric (in $p$) convex mathematical programs whose objective functions are   depending on the benefits of the sectors related.
Such a model frequently arises  in practice and has been considered by several authors, see e.g. \cite{HMD1} and \cite{Ko2} Chapter 9.
We employ the fact that under some  mild assumptions     both  maps $S$ and $D$ possess certain monotonicity properties, which allow us to apply the available solution-methods e.g. \cite{FP1} Chap12,  \cite{Ko1}  
to the variational inequality  or complementarity problems for finding  an equilibrium point of the model. However, in general, this model may have many solutions, so it is ill-posed. In order to overcome the ill-posedness, we propose a regularization algorithm based upon the bilevel optimization for finding an equilibrium price of the model that is the nearest to a given guessed or desired  equilibrium point.

The paper is organized as follows.   The next section contains  some definitions for monotonicity properties for an operator and lemmas to be used in the forthcoming sections.  In the third section  we describe the model in detail. The section fourth is devoted to describe an algorithm for minimizing the distance from a guessed or desired equilibrium point to the equilibrium point-set of the model. We close the paper with some reports on numerical results for the model. The obtained computational results on many randomly generated data show that the proposed algorithm works well for this bilevel problem.

\section{Priliminaries } 
First, let us recall some definitions, see e.g. \cite{Tuy}, pages 39, 40, and lemmas that will be used in the forthcoming sections.

Let $C \subseteq \mathbb {R}^n $ be a convex set and $f : \mathbb{R}^n \to \mathbb{R} \cup{\{+\infty\}}$ such that
$f(x) < +\infty$ for all $x\in C$. The function $f$ is said to be convex on $C$, if
$$f\big (\lambda x+(1-\lambda )y\big) \leq \lambda f(x) + (1-\lambda)f(y) \ \forall x,y \in C, \ \forall \lambda \in [0, 1].$$
$f$ is said to be strongly convex with modulus $\gamma > 0$ (shortly $\gamma-$ strongly convex) on $C$ if 

$$f\big (\lambda x+(1-\lambda )y\big) \leq \lambda f(x) + (1-\lambda)f(y) - \gamma\|x-y\|^2 \ \forall x,y \in C.$$
The function $f$ is said to be concave (resp. strongly concave) on $C$ if $- f$ is convex (resp. strongly convex) on $C$. 

The subdifferential of $f$ at a point $x$, denoted by $\partial f(x)$, is defined  as
$$\partial f(x):=\{ u \in \mathbb{R}^n : \langle u, y-x \rangle + f(x) \leq f(y) \ \forall y\}.$$
It is well known that if $f$ is differentiable at $x$, then $\partial f(x) \equiv\{\nabla f(x)\}$.

Let $ T: \mathbb{R}^n \to\mathbb{R}^n $ such that $T(x) \not=\emptyset $ for every $x\in C.$ 
The following concepts for monotonicity of an operator can be found, for example, \cite{BC1} page 293.

(i) $T$ is said to be nonexpansive on $C$ if
$$\| T(x)-T(y)\| \leq \|x-y\| \ \forall x,y \in C.$$
A typical example for nonexpansive mappings is the metric projection that mapps every $x\in \mathbb{R}^n$ onto $C$ by taking
$$P_C(x) :=\{y\in C: \|x-y\| \leq \|x-z\| \ \forall z\in C\};$$
It is well known, see e.g. \cite{BC1} page 61, that if $C$ is closed, convex, then $P_C(x)$ is singleton and nonexpansive on
$ \mathbb{R}^n $, that is
$\|P_C(x)-P_C(y)\| \leq \|x-y\|$ for every $x, y \in \mathbb{R}^n$.

(ii) $T$ is said to be monotone on $C$ if
$\langle T(x) -T(y), x-y\rangle \geq 0\ \forall x,y \in C$.

(iii) $T$ is said to be strongly monotone on $C$ with modulus $\eta > 0$ if
$\langle T(x) -T(y), x-y\rangle \geq \eta \|x-y\|^2 \ \forall x,y \in C$.

A typical example for monotone (resp. strongly monotone)  operators is the gradient of a lower semi continuous convex (resp. strongly convex) function.

(iv) $T$ is said to be co-coersive (inverse-strongly monotone ) \cite{ZM1} on $C$ with modulus $\beta >0$ if
$\langle T(x) -T(y), x-y\rangle \geq \beta\|T(x)-T(y)\|^2\ \forall x,y \in C$.

It is known \cite{ZM1} that the gradient of a convex function on an open convex set containing $C$ is inverse-strongly monotone whenever the gradient is Lipschitz on $C$.

Below are the lemmas that will be used in the proof of the convergence theorem for the algorithm to be described.

\begin{lemma}\label{M1aa} \cite{Xu1} Suppose that $\{\alpha_k\}$ is a sequence of nonnegative numbers such that 
$$\alpha_{k+1}\leq (1-\lambda_k)\alpha_k + \lambda_k \delta_k + \sigma_k \ \forall k,$$
where

(i) $0 <\lambda_k < 1, \ \sum_{k=1}^{+\infty} \lambda_k = +\infty$;

(ii) $\limsup_{k\to +\infty} \delta_k \leq 0$;

(iii) $\sum_{k=1}^{+\infty} |\sigma_k| < +\infty$.

Then $\lim \alpha_k = 0$ as $k\to +\infty$.
\end{lemma}

\begin{lemma}\label{M1a} (\cite{BC1} Theorem 4.17) Let $C$ be a nonempty closed subset in a Hilbert space $\mathcal{H}$ and 
$T : C \to \mathcal{H}$ be a nonexpansive mapping. Let $\{x^k\}$ be a sequence in $C$ such that $x^k \to x$, $x^k - T(x^k) \to u$. Then $x -T(x) = u$.

\end{lemma}

\begin{lemma}\label{M1b} \cite{MQ1} Let $\varphi$ be a strongly convex differentiable function with modulus $\beta$ and $L$-Lipschitz on $C$.
Then, for any $\alpha >0$ it holds that
$$\|\big(x-\frac{1}{\alpha} \nabla \varphi(x) \big) - \big (y-\frac{1}{\alpha} \nabla \varphi(y) \big)\|^2 \leq (1-\frac{2\beta}{\alpha} +\frac{L^2}{\alpha^2}) \|x-y\|^2 \ \forall x, y \in C.$$

\end{lemma}

\section{The model}

As we have mentioned, the model to be solved has two agents, each of them has a strategy set depending on a parameter $p$ laying in a closed convex set $P$.
Given a  vector $p\in P$, the first agent  determines his/her strategies
$S(p) \subset X \subseteq \mathbb{R}^n_+$, while the second one determines his/her strategies $D(p) \subset X $, where $X$ is a given closed convex set. Then we  define  the strategy $F$ for the model
 by taking
$F(p)= S(p)-D(p)$.
We recall that vector $p^* \in P$ is said to be an equilibrium point if it is a solution of the multivalued variational inequality problem
$$\text{Find}\ p^*\in P: q^* \in F(p^*): \langle q^*,p-p^*\rangle \geq 0\  \forall p\in P. \eqno MVI(P,F)$$
In the case $P$ is a closed convex cone, this problem is reduced to the complementarity one
$$ p^* \in P: \exists q^* \in F(p^*), \langle q^*, p^*\rangle = 0. \eqno MCP(P, F)$$
In what follows we suppose that both $S$ and $D$ are given respectively as the optimal solution-sets of the parametric (in $p$) mathematical programs

\begin{equation}\label{S1}
\max \{p^T x - c(x)  \ | \ \text{s.t}. \ x\in X \}
\end{equation}
and
\begin{equation}\label{D1}
\min \{p^T x + t(x) \ | \ \text{s.t}. \ x\in X, u(x) \geq M > 0\}, 
\end{equation}
respectively, where $c(x)$ is the cost for strategy $x$, for example the producing cost (including the environmental cost), while $u(x)$ is the utility of the second agent and $t(x)$ is the tax for  $x$.
We suppose that the functions $c$, $t$ and $-u$ are convex. The convexity of function $c$ (resp. $t$) means that the cost for producing (resp. the tax for using) a unique commodity increases as the amount gets larger. Thus Problems (\ref{S1}) and (\ref{D1}) are convex.

Note that, in \cite{Ko2}, page 152, the mapping $D$ has been defined as the optimal  solution-set of the mathematical program 
\begin{equation}\label{D2}
 \text{arg}\max \{ u(x) \ | \ \text{s.t.}\ p^Tx \leq \omega, x\in X\},
\end{equation}
where $u$ is an $\alpha$-positively homogeneous concave function and $\omega > 0 $ is the budget. It is easy to see that the solution-set of problem   (\ref{D2})
is just the solution-set of the convex mathematical program 
\begin{equation}\label{D3}
\min \{ p^Tx - \frac{\omega}{\alpha}\ln u(x) \ | \ \text{s.t.}\ x\in X\},
\end{equation}
which is a special form  of  (\ref{D1}).

\section{ A regularization algorithm by bilevel optimization}
 Consider  the  mapping \cite{Ro2}
\begin{equation}\label{D4}
	T := (I+\lambda \Phi)^{-1},
\end{equation}\label{D4}
where $\Phi := F+ N_P$,  $I$ is the identity, $N_P$ is the normal cone of $P$
and $\lambda>0$.  
Since  $F$ is monotone  and the normal cone $N_P$ is maximal monotone,  $\Phi$ is  maximal monotone. If, in addition, $F$ is single valued, the mapping $T$ is single valued nonexpansive, see  \cite{Go,Ro2}.  

In order to develop an algorithm for minimizing the distance function over the equilibrium points of the model, we consider a special case when both $c$ and $t$ are differentiable, strongly convex. Then  we can find an equilibrium of the model  by finding  a fixed point of the  projection mapping
$ P_P\big( p-\lambda F(p)\big)$, which is easy to compute when  $P$ has a simple structure such as $\mathbb{R}^n_+$ or a rectangle (often in practice). In fact, we have the following lemma.  

\begin{lemma}\label{M4a} \cite{HMD1} Suppose that $c$ and $t$ are differentiable, strongly convex with modulus $\mu_c$ and $\mu_t$ respectively on an open convex set containing $X$, then the mappings $S(\cdot)$, $-D(\cdot)$ are single valued inverse-strongly monotone with modulus $\mu_c$ and $ \mu_t$, and Lipschitz with constant $L_c= \frac{1}{\mu_c}$ and $L_t= \frac{1}{\mu_t}$ on $P$, respectively. Then the projection mapping $T(p):= P_P\big( p-\lambda F(p)\big)$ is nonexpansive for every $0\leq\lambda \leq 2\mu_F $, with $\mu_F= \frac{1}{2}\min\{\mu_c, \mu_t\}$ and its fixed point-set coincides with the solution-set of variational inequality defined as
$$ \text {find}\ p^*\in P: \langle F(p^*), p-p^*\rangle \geq 0. \ \eqno VI(P, F)$$
\end{lemma}

  Note that the fixed point-set of the nonexpansive mapping $T$, in general, is not a singleton, so the problem of finding its a fixed point is not well-posed. In order to overcome this ill-posedness, we use a bilevel optimization approach. This regularization  approach has been used by several authors, see e.g. \cite{DY1,DHM1}.
For this purpose, one can choose a strongly convex continuously differentiable function $ f$ on $\mathbb{R}^n$, for example, $f(p) := \|p -p^0\|^2$ with $p^0$ being given (plays as a guessed or desired) equilibrium price. Then we consider the following mathematical programming problem
$$\min f( p) \ \text{s.t.} \ p\in Fix(T). \eqno(BOP)$$
Since $Fix(T)$ is closed convex and $f$ is strongly convex, this problem always has a unique optimal solution, however, since the constrained set of Problem (BOP) is given implicitly as the fixed point-set of $T$, the existing algorithms for convex programming cannot be applied directly to this problem.

Now we describe an algorithm for solving problems (BOP) and study its convergence analysis. We suppose that the metric projection onto the closed convex set $P$ can be computed with a reasonable effort, for example, $P = \mathbb{R}^n_+ $ or a rectangle.\\

The following algorithm is a combination between the gradient one for minimizing the strongly convex function $f$ and the Mann-Krasnosel'skii iterative scheme for approximating a fixed point of the nonexpansive mapping $T$.

ALGORITHM 

{\it Initialization}. Choose $0\leq \eta \leq 2\mu_F$ and two sequences of positive numbers $\{\lambda_k\}$, $\{\alpha_k\}$ decreasing to zero. In addition,
$$\begin{array}{lll}
\sum_{k=1}^{+\infty} \lambda_k \alpha_k &=& +\infty;\\ 
\sum_{k=1}^{+\infty} |\lambda_{k+1} - \lambda_k| &<& +\infty;\\
\sum_{k=1}^{+\infty} |\alpha_{k+1} - \alpha_k| &< & +\infty.
\end{array}
$$
Pick $p^1\in P$ and let $k = 1.$\\

{\bf Iteration} $k = 1,2...$ Having $p^k\in P$, and $g^k = \nabla f (p^k)$, compute
$$ q^k := \text{argmin}\{ \langle g^k, y-p^k\rangle + \frac{1}{2\alpha_k}\|y-p^k\|^2\ : y \in P\} 
= P_P(p^k - \alpha_k g^k).$$
Define the next iterate
$$p^{k+1}:= \lambda_k q^k + (1-\lambda_k) T(p^k)$$
with 
$$T(p^k):= P_P\big (p^k -\eta F(p^k) \big).$$
If $p^k = q^k = p^{k+1}$, terminate: $p^k$ is a solution to (BOP), otherwise increase $k$ by one and go
to iteration $k$.
\begin{theorem} (i) If the algorithm terminates at some iteration $k$, then $p^k$ is the solution of (BOP).

(ii) If the algorithm does not terminate, then it generates an infinite sequence $\{ p^k\}$ that converges to the unique solution of (BOP).
\end{theorem}

{\bf Proof.} 
(i) According to the algorithm, $p^k= q^k$ means that $p^k = P_P\big (p^k - \alpha_k \nabla f(p^k) \big)$, which, by convexity of $f$, implies that $p^k$ is the minimizer of $f$ on $P$. In addition, if 
$$p^k =p^{k+1} = \lambda_k p^k + (1-\lambda_k) T(p^k),$$ 
then $p^k = T(p^k)$, which means that $p^k$ is a fixed point of $T$.

(ii) The proof is divided into several steps.

{\it Step 1}. We show that $\{p^k\}$ is bounded. Indeed, let $p^*$ be the unique solution of (BOP), then according to the algorithm and the nonexpansiveness of the projection map, we have

\begin{equation}\label{e1} \begin{array}{lll}
\|p^{k+1}-p^*\|&=& \|\lambda_k q^k +(1-\lambda_k)T(p^k)- p^*\| \\
&= & \|\lambda_k ( q^k - p^*)+ (1- \lambda_k)\big(T(p^k)-T(p^*)\big)\|\\
&\leq & \lambda_k \| q^k - p^*\|+(1-\lambda_k) \|T(p^k)-T(p^*) \|\\
&= &\lambda_k\|P_{P}\big{(}p^k - \alpha_kg^k\big{)}-P_{P}(p^*)\| + (1-\lambda_k) \|T(p^k)-T(p^*) \|\\
&\leq & \lambda_k \| p^k - \alpha_k g^k - p^*\|
+(1-\lambda_k) \|p^k - p^*\|.
\end{array}
\end{equation}
Let $g^* = \nabla f(p^*)$.
Since $g^k = \nabla f(p^k)$ and the strong convexity of $ f$ with modulus $\beta >0$, by Lemma \ref{M1b}, we have

\begin{equation}\label{e27}
\|p^k -p^* - \alpha_k (g^k -g^*)\|^2 \leq (1- 2\beta\alpha_k + L^2\alpha_k^2 ) \|p^k-p^*\|^2= (1-\gamma_k)^2
\| p^k- p^*\|^2,
\end{equation}
where
$$0 < \gamma_k = 1- \sqrt{1- 2\beta\alpha_k+L^2\alpha_k^2}<1.$$
Thus
\begin{equation} \label{e28}
\begin{array}{lll}
\|p^k-\alpha_k g^k-p^*\|&\leq &\|p^k-p^*-\alpha_k(g^k-g^*)\| +\alpha_k\|g^*\|\\
&\leq &(1-\gamma_k)\|p^k-p^*\| + \alpha_k\|g^*\|.
\end{array}
\end{equation}
It is easy to see that
\begin{equation}\label{e29}
\lim_{k\rightarrow \infty} \frac{\alpha_k}{\gamma_k} = \frac{1}{\beta}.
\end{equation}
\noindent From (\ref{e29}), it follows that there exists $k_0\in \mathbb{N}$ such that
$$\frac{\alpha_k}{\gamma_k}\leq\frac{2}{\beta} \ \forall k\geq k_0. $$
Combining (\ref{e1})-(\ref{e29}), for all $k\geq k_0$, we obtain
\begin{align*}
\|p^{k+1}-p^*\|&\leq \lambda_k(1-\gamma_k )\|p^k-p^*\| +\lambda_k\alpha_k\|g^*\| + (1-\lambda_k)\|p^k-p^*\|\\
&= (1-\lambda_k\gamma_k)\|p^k-p^*\| + \lambda_k \alpha_k\|g^*\|\\
&= (1- \lambda_k\gamma_k)\|p^k-p^*\| + \lambda_k\gamma_k\Big{(}\frac{\alpha_k}{\gamma_k}\Big{)}\|g^*\|\\
& \leq (1- \lambda_k\gamma_k)\|p^k-p^*\| + \lambda_k\gamma_k\frac{2\|g^*\|}{\beta} \\
&\leq \max \big{\{}\|p^k- p^*\|, \frac{2\|g^*\|}{\beta}\big{\}},
\end{align*}

\noindent from which, by induction, we obtain
$$\|p^{k+1}-p^*\|\leq \max\Big{\{}\|p^{k_0}- p^*\|, \frac{2\|g^*\|}{\beta}\Big{\}}.$$
\noindent Thus $\{p^k\}$ is bounded, and therefore $\{g^k\}$ and $\{q^k\}$,  $\{T(p^k)\}$ are bounded too.

{\it Step 2}. Now we prove that any cluster point of the sequence $\{p^k\}$ is a fixed point of $T$. Indeed, by the definition, the boundedness of $\{q^k\}$, $\{T(p^k)\}$ and $\lambda_k \to 0$ as $k\rightarrow \infty$, we have
\begin{equation}\label{e3} \begin{array}{lll}
\lim_{k\rightarrow \infty}\|p^{k+1} - T(p^k)\|&= &\lim_{k\rightarrow \infty} \lambda_k \|q^k -T(p^k)\|\\
& = & 0. 
\end{array}
\end{equation}

\noindent Now let k
$$K:= \sup_{k} \Big(\max(\{\|q^{k-1}\|+\|T(p^{k-1})\|, \|g^{k-1}\|\})\Big),$$
then $K < \infty$.

By using the same techniques as above, we have
$$ \begin{array}{lll}
\|p^{k+1} - p^k\| &=& \|\lambda_kq^k + (1-\lambda_k)T(p^k) - \lambda_{k-1}q^{k-1} -(1-\lambda_{k-1})T(p^{k-1})\|\\
&=& \|\lambda_k(q^k - q^{k-1}) + (1-\lambda_k)(T(p^k)-T(p^{k-1})) + (\lambda_{k-1} - \lambda_k)(-q^{k-1}+T(p^{k-1})\|\\
&\leq& \lambda_k\|P_{P}(p^k - \alpha_kg^k) - P_{P}(p^{k-1}-\alpha_{k-1}g^{k-1})\| + (1-\lambda_k)\|p^k-p^{k-1}\|\\
&&+ |\lambda_k - \lambda_{k-1}|(\|q^{k-1}\|+\|T(p^{k-1})\|)\\
&\leq & \lambda_k\|p^k - p^{k-1} - \alpha_k(g^k - g^{k-1})\| + \lambda_k|\alpha_k - \alpha_{k-1}|\|g^{k-1}\|\\
&& + (1-\lambda_k)\|p^k-p^{k-1}\| + K|\lambda_k - \lambda_{k-1}|\\
&\leq & (1-\lambda_k\gamma_k)\|p^k-p^{k-1}\| + K (|\lambda_k -\lambda_{k-1}|+|\alpha_k-\alpha_{k-1}|).
\end{array}
$$
Using (\ref{e29}), and 
$\sum_{k=1}^{+\infty} \lambda_k\alpha_k = +\infty$, one has $\sum_{k=1}^{+\infty} \lambda_k\gamma_k = +\infty$. Moreover, from the assumption $\sum_{k=1}^{+\infty} |\lambda_k -\lambda_{k-1}|< +\infty$ and $\sum_{k=1}^{+\infty} |\alpha_k -\alpha_{k-1}|< +\infty$, by Lemma \ref{M1aa},
we have $\|p^{k+1} - p^k\| \to 0$, which together with (\ref{e3}) implies
$\|p^k-T(p^k)\|\to 0$ as $k\rightarrow \infty$. Thus any cluster point of $\{p^k\}$ is a fixed point of $T$.

{\it Step 3}. We show that $\{p^k\}$ converges to $p^*$. In fact, by definition of $p^{k+1}$ and $q^k$, the nonexpansiveness of the projection and $T$, we can write
\begin{equation} \label{e40}
\begin{array}{lll}
\|p^{k+1} -p^*\|^2 &=& \|\lambda_kq^k + (1-\lambda_k)T(p^k) -p^*\|^2\\
&=& \|\lambda_k(q^k-p^*) +(1-\lambda_k)(T(p^k)-T(p^*))\|^2\\
&=& \|\lambda_k\Big{(}P_{P}(p^k-\alpha_kg^k)-P_{P}(p^*)\Big{)} + (1-\lambda_k)(T(p^k)-T(p^*))\|^2\\
&\leq& \lambda_k\|p^k-\alpha_kg^k-p^*\|^2 + (1-\lambda_k)\|p^k-p^*\|^2\\
&=& \lambda_k\|p^k - \alpha_kg^k - p^* + \alpha_kg^* - \alpha_kg^*\|^2 + (1-\lambda_k)\|p^k-p^*\|^2\\
&=& \lambda_k\|p^k - p^* -\alpha_k(g^k-g^*)\|^2 +\lambda_k\alpha_k^2\|g^*\|^2 \\
&&+ 2\lambda_k\alpha_k\langle -g^*, p^k - p^* -\alpha_kg^k +\alpha_kg^*\rangle + (1-\lambda_k)\|p^k-p^*\|^2\\
&\leq&\lambda_k\|p^k - p^* -\alpha_k(g^k-g^*)\|^2 + 2\lambda_k\alpha_k\langle -g^*, p^k - p^* -\alpha_kg^k \rangle \\
&& + (1-\lambda_k)\|p^k-p^*\|^2.
\end{array}
\end{equation}
From (\ref{e40}), using Lemma \ref{e27} and recall that $0<1-\gamma_k <1$, we obtain
\begin{equation}\label{e4}
\|p^{k+1} -p^*\|^2\leq (1-\lambda_k\gamma_k)\|p^k-p^*\|^2 + \lambda_k\gamma_k \Big{(}\frac{2\alpha_k}{\gamma_k}\Big{)}\langle -g^*, p^k-p^* - \alpha_kg^k\rangle.
\end{equation}
Let $\{p^{k_j}\}$ be a subsequence of $\{p^k\}$ such that
$$\limsup_{k\rightarrow \infty}\langle -g^*,p^k-p^*\rangle = \lim_{j\rightarrow \infty}\langle -g^*,p^{k_j}-p^*\rangle .$$
Using a subsequence of $\{p^{k_j}\}$, if necessary, we may assume that $p^{k_j} \to \bar{p}$ as $j \to +\infty$. Then, by Step 2, $\bar{p}$ is a fixed point of $T$. Recall that $\{g^k\}$ is bounded and $\alpha_k \to 0$, using (\ref{e29}), we obtain 

\begin{equation}\label{e41}
\limsup_{k\rightarrow \infty}\frac{2\alpha_k}{\gamma_k}\langle -g^*,p^k-p^*-\alpha_kg^k\rangle = \lim_{j\rightarrow \infty}\frac{2\alpha_{k_j}}{\gamma_{k_j}}\langle -g^*,p^{k_j}-p^*- \alpha_{k_j}g^{k_j}\rangle =\frac{2}{\beta} \langle -g^*, \bar{p} -p^*\rangle.
\end{equation}
Now, since $p^*$ is the solution of the strongly convex optimization problem: 
$$\min \{ f(p): p\in Fix(T)\}$$ 
and $g^* = \nabla f(p^*)$, we have 
$$\langle g^*, p-p^*\rangle \geq 0 \ \forall p\in Fix(T),$$ 
in particular, 
$$\langle g^*, \bar{p}-p^*\rangle \geq 0,$$ 
which together with (\ref{e41}) implies that
\begin{equation} \label{e42} 
\limsup_{k\rightarrow \infty}\frac{2\alpha_k}{\gamma_k}\langle -g^*,p^k-p^*-\alpha_kg^k\rangle \leq 0. 
\end{equation}
Thus, from (\ref{e4}) and (\ref{e42}), by applying Lemma \ref{M1aa} with
$\sigma_k \equiv 0$, we can deduce that $\|p^k -p^*\| \to 0$ as $k\rightarrow \infty$. \hfill $\Box$
\begin{remark}\label{rm51}
(i) Case (a) occurs only if the unique minimizer of $f$ over $P$ is also a fixed point of $T$ on $P$.

(ii) When $P = \mathbb{R}^n_+$, the $j$-component of vector $q^k$ is given by 
\begin{equation*}
q_j^k= \begin{cases}
p_j^k- \alpha_k g_j^k \ \ \textrm{if} \ p_j^k- \alpha_kg_j^k \geq 0,\\ 
0 \ \ \ \hspace{1.1cm} \ \textrm{if} \ \ p_j^k- \alpha_kg_j^k < 0.
\end{cases}
\end{equation*}
(iii) When $P$ is a rectangle given as 
$P:= \{p^T:= (p_1,...,p_n): \ a_j \leq p_j \leq b_j, \  \forall j\}$, then 
\begin{equation*}
	q_j^k= \begin{cases}
		p_j^k- \alpha_kg_j^k \ \ \textrm{if} \ a_j \leq p_j^k- \alpha_kg_j^k \leq b_j,\\ 
		a_j \ \ \hspace{1cm} \ \textrm{if} \ \ p_j^k- \alpha_kg_j^k < a_j,\\
		 b_j \ \ \hspace{1cm} \ \textrm{if} \ \ p_j^k- \alpha_kg_j^k > b_j.
	\end{cases}
\end{equation*}

\end{remark}

\section{Computational results} In this section we consider the model for two cases. For the first case we tale  $P := \mathbb{R}^n_+$, while for the second one $P $ is a rectangle. For these cases the projections onto them have the closed form. In both cases we take  
$c(x) = x^T C x , \ t(x) = x^TB x, u(x) = l^Tx$ and
$f(p) := \|p-p^0\|^2$, where all input data
$p^0$  is randomly generated in the interval $[0, 100]$. We take 
$C:= C_1^T C_1, B :=B_1^T B_1$ whose entries are randomly generated in the interval $[-10, 10]$. Clearly these matrices are symmetric and positive definite.
The feasible domain is given as $X:= \{ x \geq 0: Ax \leq b\}$ where $A$ is a $m\times n$ matrix, $b \in \mathbb{R}^m$ whose all entries are randomly chosen in $(0, 20)$.\

For this model we choose the parameter $\lambda_k = \alpha_k = \frac{1}{\sqrt{k+1}}$ for all $k$.

We tested the model with Python 3.10
on a computer with the processor: AMD Ryzen 5 1600 Six- Core Processor 3.20 GHz with the installed memory (RAM): 16.0 GB. \\

We stopped the computation when $ \frac{\|p^{k+1}-p^k\|}{max{\{\|p^{k+1}\|, 1\}}} < \epsilon$, with $ \epsilon = 10^{-4}$;\\


The number of the tested problems is 10. In the tables below we use the following headings:
\begin{itemize}

\item {Average Times: the average CPU-computational times (in second);}

\item {Average Iteration: the average number of iterations.}

\end{itemize}
\begin{center}
Table 1. Computed Results for  $P := \mathbb{R}^n_+$

\begin{tabular}{|c|c|c|c|}
\hline

\ \ \hspace{0.25cm} n \hspace{0.5cm} & \hspace{0.5cm} m \hspace{0.5cm} & \hspace{0.25cm}Average Times\hspace{0.25cm}
& \hspace{0.25cm} Average Iteration\ \ \ \ \\
\hline
5&3& 1.777 & 87.6 \ \ \ \ \\
\hline
10&8& 2.479 & 107.1 \ \ \ \ \\
\hline 
30&20& 4.043 & 127.4 \ \ \ \ \\
\hline 

50&30& 7.391 & 139.7 \ \ \ \ \\
\hline 
100&80& 11.341 & 152.4 \ \ \ \ \\
\hline 
\end{tabular}
\end{center}

\begin{center}

Table 2. Computed Results when $P$ is a rectangle

\begin{tabular}{|c|c|c|c|}
\hline

\ \ \hspace{0.25cm} n \hspace{0.5cm} & \hspace{0.5cm} m \hspace{0.5cm} & \hspace{0.25cm}Average Times\hspace{0.25cm}
& \hspace{0.25cm} Average Iteration\ \ \ \ \\
\hline
5&3& 1.632 & 123.3 \ \ \ \ \\
\hline
10&8& 2.223 & 152.9 \ \ \ \ \\
\hline 
30&20& 5.026 & 176.5 \ \ \ \ \\
\hline 

50&30& 7.623 & 195.3 \ \ \ \ \\
\hline 
100&80& 9.651 & 236.5 \ \ \ \ \\
\hline

\end{tabular}
\end{center}

\textbullet \ The graph of $\frac{\|x^{k+1}-x^k\|}{\max\{\|x^k\|,1\}}$ with respect to $Iteration$ when $P := \mathbb{R}^n_+$.

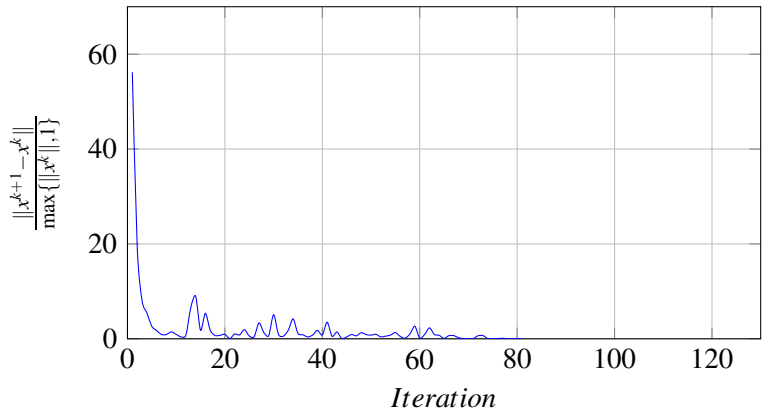
\begin{center}

\begin{tikzpicture}
\begin{axis}[
    xlabel={$Iteration$},
    ylabel={$\frac{\|x^{k+1}-x^k\|}{\max\{\|x^k\|,1\}}$},
    xmin=0, xmax=130,
    ymin=0, ymax=70,
    grid=major,
    width=10cm,
    height=6cm,
    cycle list name=exotic,
    legend pos=north west,
    smooth,
    mark=*,
    ]
\addplot[blue] coordinates {
(1, 56.197576595493783)
(2, 19.386802051356284)
(3, 8.0577035249053682)
(4, 5.4490687794272224)
(5, 2.6161034742566844)
(6, 1.7446395659057385)
(7, 0.8866027093423826)
(8, 0.9033917981728566)
(9, 1.4626463995742842)
(10, 0.945140374250422)
(11, 0.411258439264986)
(12, 0.785247844917183)
(13, 6.581924537501174)
(14, 8.986469116439993)
(15, 1.796417613651980)
(16, 5.377955463303499)
(17, 1.775412128435022)
(18, 0.682858756998600)
(19, 0.737581797610187)
(20, 0.954743237313193)
(21, 0.032882479997257)
(22, 1.008557193512305)
(23, 0.816375016093845)
(24, 1.919917717199210)
(25, 0.590436042560398)
(26, 0.420821261985075)
(27, 3.350381423856556)
(28, 1.321763802327592)
(29, 0.660826297934315)
(30, 5.073230452768991)
(31, 0.931219190098964)
(32, 0.551420926021653)
(33, 1.778818600027989)
(34, 4.211189036050485)
(35, 1.123286447899532)
(36, 0.845352805207263)
(37, 0.393499661270131)
(38, 0.776533487023274)
(39, 1.772878160953895)
(40, 0.794434791097754)
(41, 3.504974988182403)
(42, 0.523402178017264)
(43, 1.448570285656420)
(44, 0.055845688592988)
(45, 0.397652037116310)
(46, 0.893634622915741)
(47, 0.625927437871665)
(48, 1.286248401496956)
(49, 0.880599264210281)
(50, 0.770113345364572)
(51, 0.949992577534288)
(52, 0.402615093625038)
(53, 0.563188272665339)
(54, 0.778341135187796)
(55, 1.309183544930329)
(56, 0.507005782830285)
(57, 0.125355631734817)
(58, 1.025427839420118)
(59, 2.654634366930886)
(60, 0.126925700846468)
(61, 0.896625547444974)
(62, 2.306512983352018)
(63, 0.916435842990113)
(64, 0.724916626359238)
(65, 0.085142150521579)
(66, 0.659906826026594)
(67, 0.674488306429247)
(68, 0.209050302574725)
(69, 0.060230072104878)
(70, 0.005979657313182)
(71, 0.064166915850879)
(72, 0.652131715579664)
(73, 0.674488306429277)
(74, 0.033124134939245)
(75, 0.060230072104878)
(76, 0.001618606862525)
(76, 0.001473001255056)
(77, 0.123510690277064)
(78, 0.001092587562385)
(79, 0.005064493368685)
(80, 0.042979333278885)
(81,0.0000033157445330) 
};
\end{axis}
\end{tikzpicture}


\end{center}

\textbullet \ The graph of $\frac{\|x^{k+1}-x^k\|}{\max\{\|x^k\|,1\}}$ with respect to $Iteration$ when $P$ is a rectangle.

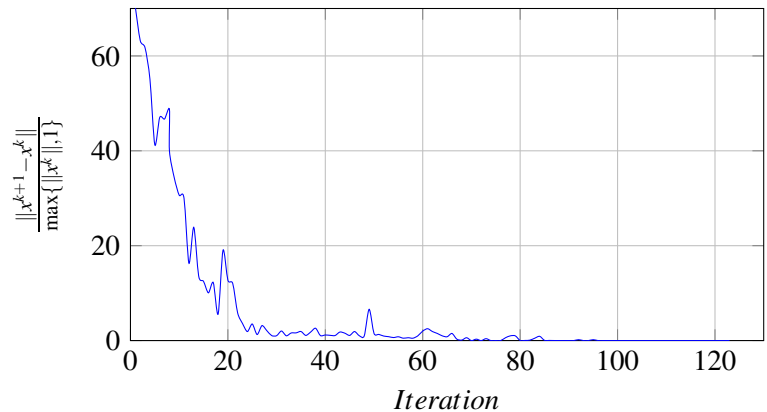
\begin{center}
    
\begin{tikzpicture}
\begin{axis}[
    xlabel={$Iteration$},
    ylabel={$\frac{\|x^{k+1}-x^k\|}{\max\{\|x^k\|,1\}}$},
    xmin=0, xmax=130,
    ymin=0, ymax=70,
    grid=major,
    width=10cm,
    height=6cm,
    cycle list name=exotic,
    legend pos=north west,
    smooth,
    mark=*,
    ]
\addplot[blue] coordinates {
            (1, 70.361503962432374)
            (2, 63.271685027466404)
            (3, 61.7040628639584647)    
            (4, 55.091775113364633)
            (5, 41.283055824137822)
            (6, 47.013515124102725)
            (7, 46.709960789382175)
            (8, 48.750902128678029)
            (8, 40.018777420731457)
            (9, 34.2843812397948774)
            (10, 30.657490801902553)
            (11, 29.966992827853034)
            (12, 16.325523186266937)
            (13, 23.887854696644959)
            (14, 13.628259059495304)
            (15, 12.485170588150334)
            (16, 10.0601870806061826)
            (17, 12.2228141074484144)
            (18, 5.59615214491199433)
            (19, 19.0332137830750387)
            (20, 12.7366629469471533)
            (21, 12.0145009537188543)
            (22, 5.84926353033850578)
            (23, 3.60544837767435413)
            (24, 1.89817625116920845)
            (25, 3.49939645665569593)
            (26, 1.27617537076919785)
            (27, 3.15200916495021083)
            (28, 2.08332437804255766)
            (29, 1.04551800018709782)
            (30, 1.02480144098246169)
            (31, 2.01349725358281846)
            (32, 1.00735568272943477)
            (33, 1.60403723991614766)
            (34, 1.62260947888172785)
            (35, 1.91328289965527195)
            (36, 1.08560848907052355)
            (37, 1.80045546124899757)
            (38, 2.60065947872556932)
            (39, 1.05121239831482166)
            (40, 1.17644516199993866)
            (41, 1.10080460771859812)
            (42, 1.05929137721377999)
            (43, 1.80082780732787164)
            (44, 1.56070545008971444)
            (45, 1.08379512781262442)
            (46, 1.90579645585398775)
            (47, 1.00956317204050872)
            (48, 1.00067414639504867)
            (49, 6.60077620455534496)
            (50, 1.50814240308613963)
            (51, 1.30019625350127718)
            (52, 0.95002731569965771)
            (53, 0.82099298376583482)
            (54, 0.63618703515355706)
            (55, 0.80097455106135548)
            (56, 0.50057335908819033)
            (57, 0.60036357268610525)
            (58, 0.50044470106282922)
            (59, 1.00088286307138653)
            (60, 2.00096024478451782)
            (61, 2.50033379737210236)
            (62, 1.90037308318832220)
            (63, 1.50023069855500834)
            (64, 1.00026844215742016)
            (65, 0.80026780956759947)
            (66, 1.50005391149374627)
            (67, 0.30091554858260727)
            (68, 0.10093825378704833)
            (69, 0.60026849044646848)
            (70, 0.00022586407988063)
            (71, 0.30049672067680706)
            (72, 0.00021414656110984)
            (73, 0.40013392886567897)
            (74, 0.00012698642823281)
            (75, 0.00012057504309517)
            (76, 0.00011464059789272)
            (77, 0.60109136193347992)
            (78, 1.00010402081549434)
            (79, 1.00019925836485555)
            (80, 0.02064338052384257)
            (81, 0.03016191344338831)
            (82, 0.09019708640373678)
            (83, 0.50068018721530116)
            (84, 0.90000093793116344)
            (85, 0.05013505521480649)
            (86, 0.03017936634880985)
            (87, 0.01012877487425244)
            (88, 0.00511982994660081)
            (89, 0.00003119618007719)
            (90, 0.00013470983171338)
            (91, 0.00015507681596935)
            (92, 0.19012971446802801)
            (93, 0.00012632990495582)
            (94, 0.00012383189550568)
            (95, 0.19036933215888776)
            (96, 0.00035621410328289)
            (97, 0.00043783638183373)
            (98, 0.00013319934733507)
            (99, 0.00032080032055373)
            (100, 0.0003101644935658)
            (101, 0.0001071483734812)
            (102, 0.0001030004955378)
            (103, 0.0000791904219966)
            (104, 0.0001812509702196)
            (105, 0.0001072508030569)
            (106, 0.0001049151141014)
            (107, 0.0001485948651565)
            (108, 0.0001008971020315)
	    (109, 0.0012351069027709)
            (110, 0.0012094829280184)
            (111, 0.0011846485032277)
            (112, 0.0011605715064289)
            (113, 0.0011372214340795)
            (114, 0.0011145693034212)
            (115, 0.0010925875623859)
            (116, 0.0010712500058902)
            (117, 0.0010505316969511)
            (118, 0.0010304088941866)
            (119, 0.0001070147851832)
            (120, 0.0002696804383054)
            (121, 0.0004999907006612)
            (122, 0.0001008971020315)
            (123, 6.1346203403336143e-05)
        };
\end{axis}
\end{tikzpicture}
\end{center}

The computational results reported in the above tables show that  the proposed algorithm works well for this class of  equilibrium models.

{\bf Conclusion remark} We have considered a class of price equilibrium models, where there are two sectors whose strategy-sets are given implicitly as the solution-sets of certain parametric optimization problems. We have used the fact that the equilibrium points coincide with the fixed points of a nonexpensive mapping. In order to avoid the ill-posedness of the model we have proposed an algorithm for optimizing the distance from a given point to the fixed point-set. Computational results computed from many randomly data show that the algorithm works well for this model. 

{\bf Author Contributions} The contributions of the authors for the article as:
The first named author N. N. Hai contributed the proof of the convergence results for  the algorithm
The second named author  L. D. Muu suggested the subject and the original draft preparation
The authors N.V. Quy and  N. N. Hai provided  the code for  computational results.
All authors contributed to writing–review and editing  supervision.

{\bf Competing interests} The authors declare that they have no conflict of interest.

\end{document}